\documentclass{article}
\usepackage{verbatim}
\usepackage{moreverb}
\usepackage{url}
\usepackage{amsmath, amssymb, setspace, mathptmx}
\usepackage{color}
\usepackage{appendix}
\author{Vadim Zaliva, lord@crocodile.org}
\date{July 17, 2012}
\title{Constructing an orthonormal set of eigenvectors for DFT matrix using Gramians and determinants}
\usepackage{graphicx}
\usepackage{listings}
\usepackage{rotating}
\usepackage[colorlinks=false,bookmarks=true,pdfauthor={Vadim Zaliva lord@crocodile.org},
            pdftitle={Constructing an orthonormal set of eigenvectors for DFT matrix using Gramians and determinants},
            pdftex]{hyperref}

\newcommand{\mathsym}[1]{{}}
\newcommand{\unicode}[1]{{}}
\newcommand{\norm}[1]{\lVert#1\rVert}
\newcommand{\ip}[2]{\left\langle{#1},{#2}\right\rangle}

\begin{document}
\lstset{language=XML,basicstyle=\small,markfirstintag=true,numbers=left,numbersep=5pt}

\maketitle

\begin{abstract}
  The problem of constructing an orthogonal set of eigenvectors for a
  DFT matrix is well studied. An elegant solution is mentioned by
  Matveev in \cite{matveev2001}. In this paper, we present a distilled
  form of his solution including some steps unexplained in his paper,
  along with correction of typos and errors using more consistent
  notation. Then we compare the computational complexity of his
  method with the more traditional method involving direct
  application of the Gram-Schmidt process. Finally, we present our
  implementation of Matveev's method as a Mathematica module.
\end{abstract}

\section{Definitions}

The normalized matrix for discrete Fourier transform (DFT) of size $n$
is defined as:

\begin{equation}
\Phi_{jk}(n)=\frac{1}{\sqrt{n}}w^{jk}, \quad j,k=0,\ldots,n-1, \quad
w=e^{i\frac{2\pi}{n}}
\end{equation}

In some literature, an alternative definition is used where
$w=e^{-i\frac{2\pi}{n}}$. It should be possible to adopt the algorithm
described here with some minimal changes.

Throughout this paper, unless explicitly specified, we will use 0-based
indices for matrices and arrays.

The scaling factor $\frac{1}{\sqrt{n}}$ ensures that $\Phi$ is
unitary. An important for us property of a unitary matrix is that its
eigenvectors corresponding to different eigenvalues are
orthogonal\cite{strang2006linear}.\footnote{There is a typo in
  \cite{strang2006linear} stating that they are ``orthonormal'', while
  it should read ``orthogonal'' or ``can be chosen orthonormal''
  instead. It has been reported to the author and acknowledged by
  him.}

If $e_k$ is an eigenvector of $\Phi$ with associated eigenvalue $\lambda_k$ then
by the definition of an eigenvector $\Phi e_k=\lambda_k e_k$. It is
also a property of eigenvectors that $\Phi^qe_k=\lambda^q_k e_k$. In
\cite{candan2011} it has been shown that $\Phi^4=I$. This gives us $I
e_k = \lambda^4 e_k$. From that, it follows that $\lambda^4_k=1$, and
eigenvalues of DFT matrix are fourth roots of unity:

\begin{equation}
\label{eq:lambda}
\lambda=(1,i,-1,-i)
\end{equation}

The well known \cite{matveev2001,candan2011} spectral decomposition of
$\Phi$ into four orthogonal projections can be defined as:

\begin{equation}
\label{eq:pk}
p_k=\frac{1}{4}\sum_{j=0}^{3}{(-i)^{jk}\Phi^j},\quad k=0,\ldots,n-1
\end{equation}

Each projection matrix corresponds to one of four the possible
eigenvalues from equation~\eqref{eq:lambda}.  The columns of each
projection matrix are eigenvectors of $\Phi$ sharing the same
eigenvalue.  

As shown in \cite{candan2011}, the multiplicity of the eigenvalue with a
value of $\lambda_k$ is equal to the trace of $p_k$. However, we can use 
simpler formulae from \cite{matveev2001} to calculate the multiplicity
of $m_k$:

\begin{equation}
\label{eq:mult}
\begin{split}
m_o &= \lfloor\frac{n+1}{4}\rfloor,\quad \text{associated with } \lambda_0=1 \\
m_1 &= \lfloor\frac{n+2}{4}\rfloor,\quad \text{associated with } \lambda_1=i \\
m_2 &= \lfloor\frac{n+3}{4}\rfloor-1,\quad \text{associated with } \lambda_2=-1 \\
m_3 &= \lfloor\frac{n}{4}\rfloor+1,\quad \text{associated with } \lambda_3=-i \\
\end{split}
\end{equation}

where $\lfloor \ldots \rfloor$ in the equation above denotes the
\textit{floor} function. Note that for convenience, $\lambda_k$ is
defined so that $\lambda_k=i^k$.

Finally, following \cite{matveev2001}, we define $v(m,k)$, where 
$m,k=0,\ldots,n-1$ as an $n$-dimensional vector which is equal to the
$m$-th row (or $m$-th column due to matrix symmetry) of $p_k$:

\begin{equation}
\label{eq:v}
v(m,k) = ([p_k]_{0,m},[p_k]_{1,m},\ldots,[p_k]_{n-1,m})
\end{equation}

In the formula above, $[p_k]_{m,n}$ denotes an element at row $m$ and column $n$ of
projection matrix $p_k$. The $p_k$ per equation~\eqref{eq:pk} could be expanded as:

\begin{equation*}
p_k=\frac{I+(-i)^k\Phi+(-i)^{2k}\Phi^2+(-i)^{3k}\Phi^3}{4}
\end{equation*}

This allows us to write a formula, computing an element of $p_k$ at
position $(j,m)$.

\begin{equation}
\label{eq:vj}
[p_k]_{j,m}=\frac{
\delta_{j,m}+
(-i)^k\frac{w^{jm}}{\sqrt(n)}+
(-1)^k\delta_{(j+m\bmod n),0}+
(-i)^{3k}\frac{w^{-jm}}{\sqrt(n)}}
{4}
\end{equation}

Using this, we can express the $j$-th element of a vector $v(m,k)$ from
equation~\eqref{eq:v} as $v_j(m,k)=[p_k]_{j,m}$.

It should be noted that equation~\eqref{eq:vj} differs slightly from
the similar equation (22) in \cite{matveev2001} accounting for a
correction. The difference is in the arguments of the second Kronecker
delta, representing $\Phi^2$ which has the following form:

\begin{equation*}
\Phi^2=\left( \begin{array}{ccccc}
1 & 0 & \cdots & 0 & 0 \\
0 & 0 & \cdots & 0 & 1 \\
0 & 0 & \cdots & 1 & 0 \\
\vdots & \vdots & \cdots & 0 & 0 \\
0 & 1 & 0 & \cdots & 0 \\
\end{array} \right)
\end{equation*}

According to Matveev's formula, which incorrectly uses
$\delta_{n-j,m}$, we get the incorrect result:

\begin{equation*}
\delta_{n-j,m}=\left( \begin{array}{ccccc}
\mathbf{0} & 0 & \cdots & 0 & 0 \\
0 & 0 & \cdots & 0 & 1 \\
0 & 0 & \cdots & 1 & 0 \\
\vdots & \vdots & \cdots & 0 & 0 \\
0 & 1 & 0 & \cdots & 0 \\
\end{array} \right)
\end{equation*}

This is correct for all elements except the one at $(0,0)$, which should
be 1 instead of 0. The expression $\delta_{(j+m\bmod n),0}$, which we are
using instead, gives us the correct representation of $\Phi^2$. An
expression similar to ours is also used in \cite{candan2011}.

\section{Finding a complete system eigenvectors of $\Phi(n)$}

Each projection matrix corresponds to one of four the possible
eigenvalues from equation~\eqref{eq:lambda}.  The columns of each
projection matrix are eigenvectors of $\Phi$ sharing the same
eigenvalue.  A complete set of eigenvectors of $\Phi$ spanning $C^n$
could be constructed from columns of orthogonalized projection
matrices taking the first $m_k$ non-zero columns of $p_k$. The first
column of $p_1$ and $p_3$ are formed by all zeros and have to be
skipped. Thus, our set of $n$ eigenvectors would consist of:

\begin{equation}
\label{eq:vbasis}
\begin{split}
v(m,0), &\quad m = 0,1, \ldots m_0-1,\quad \text{associated with } \lambda_0=1\\
v(m,1), &\quad m = 1,2, \ldots m_1,\quad \text{associated with } \lambda_1=i\\
v(m,2), &\quad m = 0,1, \ldots m_2-1,\quad \text{associated with } \lambda_2=-1\\
v(m,3), &\quad m = 1,2, \ldots m_3,\quad \text{associated with } \lambda_3=-i\\
\end{split}
\end{equation}

with $m_0+m_1+m_2+m_3=n$. The proof of this using Chebychev sets could
be found in \cite{parks}.

\section{Orthonormalization}

Eigenvectors corresponding to different eigenvalues are orthogonal.
However, eigenvectors within the same projection matrix are not
guaranteed to be orthogonal, so the associated set of eigenvectors
does not possess the orthogonality property either.

A straightforward approach to get orthonormal eigenvectors as
suggested in Candan\cite{candan2011} is to apply Gram-Schmidt process
to all columns of each projection matrix.  Each projection matrix
$p_k$ will have rank $m_k$ and thus after normalization, the resulting
orthonormalized vector set will contain exactly $m_k$ non-zero
vectors.

Matveev in \cite{matveev2001} presents another approach to
constructing an orthonormal basis based on the same principles as the
Gram-Schmidt process but involving the use of Gramian matrices and
determinants.

The calculations of the orthogonal basis of $p_k$ involve $m_k$
columns of $p_k$ taken according to equation~\eqref{eq:vbasis}. We
have two cases: one for odd values of $k=1,3$ and one for even values
of $k=0,2$. Let us consider the case of even values first.

We can find a sequence of orthogonal vectors
$(e_0(k),e_1(k),\ldots,e_{m_k-1}(k))$ spanning eigenspace $p_k$ using
Gramian matrices and determinants\cite{gantmakher1959}:

\begin{equation}\label{eq:gmv}
\begin{split}
&e_0(k) = v(0,k), \\
&\ldots\\
&e_j(k)=\left|\begin{array}{cccc}
\ip{v(0,k)}{v(0,k)} & \cdots & \ip{v(0,k)}{v(j-1,k)} & v(0,k)\\
\vdots & \vdots & \vdots & \vdots \\
\ip{v(j-1,k)}{v(0,k)} & \cdots & \ip{v(j-1,k)}{v(j-1,k)} & v(j-1,k) \\
\ip{v(j,k)}{v(0,k)} & \cdots & \ip{v(j,k)}{v(j-1,k)} & v(j,k)\\
\end{array}\right|
\end{split}
\end{equation}

In the equation above, $\ip{v}{u}$ denotes the inner product of the
vectors $v$ and $u$. The determinant notation assumes generic
determinant formulation which is defined for matrices containing mixed
scalar and vector entries. The determinant could be calculated using
Laplace (cofactor) expansion.

It has been observed in \cite{matveev2001} that $p_k$ is in fact a
Gramian matrix of a set of vectors $v(m,k), m=0,\ldots,n-1$, such as
$[p_k]_{j,m}=\ip{v(j,k)}{v(m,k)}$ Using this fact, we can replace
$(j+1) \times j$ upper entries of the matrix in equation~\eqref{eq:gmv}
with corresponding entries from $p_k$:

\begin{equation}\label{eq:epk02}
\begin{split}
&e_0(k) = v(0,k), \\
&\ldots\\
&e_j(k)=\left|\begin{array}{cccc}
[p_k]_{0,0} & \cdots & [p_k]_{0,j-1} & v(0,k)\\
\vdots & \vdots & \vdots & \vdots \\
{}[p_k]_{j-1,0} & \cdots & [p_k]_{j-1,j-1} & v(j-1,k) \\
{}[p_k]_{j,0} & \cdots & [p_k]_{j,j-1} & v(j,k)\\
\end{array}\right|
\end{split}
\end{equation}

The resulting system of vectors $e_k$ is orthogonal but not yet
orthonormal. Each vector is normalized by dividing by its norm. As
shown in \cite{gantmakher1959}, the norm can be calculated by:

\begin{equation}\label{eq:ne02}
\norm{e_i}=\sqrt{G_jG_{j+1}}
\end{equation}

where $G_j,G_{j+1}$ are principal minors of $p_k$ of respective
orders. They represent Gram determinants.

For $k=1,3$, we need to take into account the fact that the first row
and the first column of $p_k,\quad k=1,3$ contain all
zeros. Therefore, for these values of $k$, the
equation~\eqref{eq:epk02} will become:

\begin{equation}\label{eq:epk13}
\begin{split}
&e_0(k) = v(1,k), \\
&\ldots\\
&e_j(k)=\left|\begin{array}{cccc}
[p_k]_{1,1} & \cdots & [p_k]_{1,j} & v(1,k)\\
\vdots & \vdots & \vdots & \vdots \\
{}[p_k]_{j,1} & \cdots & [p_k]_{j,j} & v(j,k) \\
{}[p_k]_{j+1,1} & \cdots & [p_k]_{j+1,j} & v(j+1,k)\\
\end{array}\right|
\end{split}
\end{equation}

and equation~\eqref{eq:ne02} will become:

\begin{equation}\label{eq:ne13}
\norm{e_i}=\sqrt{
\left|\begin{array}{ccc}
{}[p_k]_{1,1} & \cdots & [p_k]_{1,j}\\
\vdots & \vdots & \vdots  \\
{}[p_k]_{j,1} & \cdots & [p_k]_{j,j}\\
\end{array}\right|
\left|\begin{array}{ccc}
{}[p_k]_{1,1} & \cdots & [p_k]_{1,j+1}\\
\vdots & \vdots & \vdots  \\
{}[p_k]_{j+1,1} & \cdots & [p_k]_{j+1,j+1}\\
\end{array}\right|
}
\end{equation}

\section{Computational complexity}

Because computational complexity of Matveev's algorithm is very high,
it is not very practical for large $n$. The complexity is mostly attributed to
multiple cofactor expansion operations which have complexity of
$O(n!)$.

For comparison: obtaining a set of non-orthogonal eigenvectors and
orthonormalizing the set using the Modified Gram-Schmidt process would
take just $2n^3$ floating point operations
(FLOPS)\cite{golub1996matrix} which translates to $O(n^3)$ in Big-O
notation.

According to the equation \eqref{eq:mult} for any $n$, the
multiplicities of different projections could differ at most by
2. That means for a reasonably big $n$, the dimensionality of each of
the four eigenspaces of $\Phi$ is at approximately $\frac{n}{4}$.

Using this observation, the performance could be improved further by
a factor of four by applying the Gram-Schmidt process to $m_k\approx\frac{n}{4}$
vectors from each projection, which gives us a total cost of
$\frac{n^3}{4}$ FLOPS to orthogonalize the complete set of
eigenvectors. Although an improvement, this is still $O(n^3)$.

\section{Mathematica implementation}

Listing~\ref{mathcode} presents the full source code of the Mathematica module
constructing a complete orthonormal set of eigenvectors of DFT matrix
$\Phi(n)$. It performs all computations and returns the results in
symbolic form. The code intention was to illustrate and validate the
algorithm, and clarity and expressiveness were chosen over
performance. It has been developed and tested with Mathematica version
8.

Since Mathematica's \emph{Det} function does not work with matrices
containing both scalars and vectors, we have implemented our own
function \emph{lDet} which finds the determinant of any matrix using
cofactor expansion across rows. For the same reason, we must use our
own function \emph{rowMinor} instead of Mathematica's \emph{Minors}.

\begin{lstlisting}[caption={Mathematica module source code},
label=mathcode]
(* ::Package:: *)

BeginPackage["dfteigh`"];

dftEigen::usage = "Orthonormal basis of DFT matrix of rank N";
lDet::usage = "Determinant using Laplace expansion";
rowMinor::usage = "Row minor";

Begin["`Private`"]

Clear[rowMinor]
rowMinor[m_,r_]:=Map[Rest,Delete[m,r]]

Clear[lDet]
lDet[m_] := Module[{rows},
  rows = Length[m];
  If[rows == 1, m[[1, 1]],
   Sum[(-1)^(1+r)*m[[r, 1]]*lDet[rowMinor[m, r]], 
       {r,1,rows}]]
]

Clear[dftEigen]
dftEigen[n_] := Module[{multiplicities, vj, e},

  multiplicities = {Floor[n/4]+1, Floor[(n+1)/4], 
    Floor[(n+2)/4], Floor[(n+3)/4] - 1};

  vj[k_, j_, m_] := (1/4)*(
    KroneckerDelta[j, m] + 
    (-1)^k KroneckerDelta[Mod[j + m, n], 0] + 
    (-I)^k*Exp[(2*Pi*I*j*m)/n]/Sqrt[n]+
    (-I)^(3*k) Exp[(2*Pi*I*(-j)*m)/n]/Sqrt[n]
  );

  e[j_,k_] := Module[{g,gv,z,xn,d0,d1},
    z= If[OddQ[k],1,0];
    If[j==0, Table[vj[k,y,z],{y,0,n-1}]/Sqrt[vj[k,0+z,0+z]],
      g=Table[vj[k,y,x],{x,z, j+z},{y,z,j+z-1}];
      gv=Table[{Table[vj[k,y,x],{y,0,n-1}]},{x,z,j+z}];
      d0=Det[Table[vj[k,y,x],{x,z,j+z},{y,z,j+z}]];
      d1=Det[Table[vj[k,y,x],{x,z,j+z-1},{y,z,j+z-1}]];
      lDet[Join[g,gv,2]]/Sqrt[d0*d1]
    ]
  ];

  Map[e[#[[1]],#[[2]]] &,
  Flatten[Table[Transpose[{Range[0,multiplicities[[m+1]]-1],
  ConstantArray[m,multiplicities[[m+1]]]}],{m,0,3}],1]]
]

End[];
EndPackage[];
\end{lstlisting}

\section{Examples}

Using the Mathematica code above, we can calculate a set of eigenvectors for
a DFT matrix of order 6. Combining them as columns of matrix $O$ gives
the following matrix, approximated numerically:

\begin{equation*}
O_6=\left(
\begin{array}{cccccc}
 0.8391 & 0. & 0. & 0.5439 & 0. & 0. \\
 0.2433 & 0.5412 & 0.6533 & -0.3753 & 0.0843 & 0.2706 \\
 0.2433 & -0.2979 & 0.2706 & -0.3753 & -0.4596 & -0.6533 \\
 0.2433 & -0.4865 & 0. & -0.3753 & 0.7505 & 0. \\
 0.2433 & -0.2979 & -0.2706 & -0.3753 & -0.4596 & 0.6533 \\
 0.2433 & 0.5412 & -0.6533 & -0.3753 & 0.0843 & -0.2706 \\
\end{array}
\right)
\end{equation*}

We can verify that it diagonalizes $\Phi(6)$ by calculating:

\begin{equation*}
O_6^{-1}\Phi(6)O_6=\left(
\begin{array}{cccccc}
 1 & 0 & 0 & 0 & 0 & 0 \\
 0 & 1 & 0 & 0 & 0 & 0 \\
 0 & 0 & i & 0 & 0 & 0 \\
 0 & 0 & 0 & -1 & 0 & 0 \\
 0 & 0 & 0 & 0 & -1 & 0 \\
 0 & 0 & 0 & 0 & 0 & -i \\
\end{array}
\right)
\end{equation*}

The diagonal contains eigenvalues repeating consistently with the
associated multiplicities $m=(2,1,2,1)$ and dimensions of the
eigenspaces.

Taking the outer product of all columns of $O_6$ we can confirm that the set
is indeed orthonormal:

\begin{equation*}
\left(
\begin{array}{cccccc}
 1. & 0. & 0. & 0. & 0. & 0. \\
 0. & 1. & 0. & 0. & 0. & 0. \\
 0. & 0. & 1. & 0. & 0. & 0. \\
 0. & 0. & 0. & 1. & 0. & 0. \\
 0. & 0. & 0. & 0. & 1. & 0. \\
 0. & 0. & 0. & 0. & 0. & 1. \\
\end{array}
\right)
\end{equation*}

The Mathematica notebook, used to make calculations in the example above was:

\begin{doublespace}
\noindent\(\text{}\\
\text{Needs}[\text{{``}dfteigh$\grave{ }${''}},\text{ToFileName}[\text{NotebookDirectory}[],\text{{``}dfteigh.m{''}}]];\\
n=6;\\
\Phi =\text{Table}\left[\frac{1}{\sqrt{n}}\text{Exp}\left[\frac{2\pi *i*k*m}{n}\right],\{k,0,n-1\},\{m,0,n-1\}\right];\\
\text{eall}=\text{dftEigen}[n];\\
o=N[\text{Transpose}[\text{eall}]];\text{MatrixForm}[\text{Round}[o,0.0001]];\\
\text{MatrixForm}[\text{Round}[o,0.0001]]\\
\text{MatrixForm}\left[\text{Round}\left[\text{Inverse}[o].N[\Phi ] . o,10^{-10}\right]\right]\\
\text{MatrixForm}[N[\text{FullSimplify}[\text{Outer}[\text{Dot},\text{eall},\text{eall},1]]]]\)
\end{doublespace}

\section{Acknowledgements}

I would like to thank Lester F. Ludwig from the New Renaissance Institute
who suggested this problem to me and who has also provided guidance,
inspiration, and suggestions.

\nocite{*}
\bibliography{dfteigh}

\begin{thebibliography}{1}

\bibitem{matveev2001}
Vladimir~B Matveev.
\newblock Interwining relations between the fourier transfom and discrete
  fourier transform, the related functional identities and beyond.
\newblock {\em Inverse Problems}, 17:633, 2001.

\bibitem{strang2006linear}
G.~Strang.
\newblock {\em Linear Algebra and Its Applications}.
\newblock Thomson, Brooks/Cole, fourth edition, 2006.

\bibitem{candan2011}
{\c{C}}.~Candan.
\newblock On the eigenstructure of {DFT} matrices [{DSP} education].
\newblock {\em Signal Processing Magazine, IEEE}, 28(2):105--108, 2011.

\bibitem{parks}
J~McClellan and T~Parks.
\newblock Eigenvalue and eigenvector decomposition of the discrete fourier
  transform, 1972.

\bibitem{gantmakher1959}
F~R Gantmakher.
\newblock {\em The Theory of Matrices}.
\newblock Chelsea Publishing Co., 1959.

\bibitem{golub1996matrix}
G.H. Golub and C.F. Van~Loan.
\newblock {\em Matrix computations}, volume~3.
\newblock Johns Hopkins Univ Pr, 1996.

\end{thebibliography}
\bibliographystyle{unsrt}

\end{document}